\documentclass[11pt]{amsart}

\usepackage[a4paper,hmargin=3cm,vmargin=3cm]{geometry}
\usepackage{amsfonts,amssymb,amscd,amstext,hyperref}
\usepackage{graphicx}
\usepackage{color}
\usepackage[dvips]{epsfig}
\usepackage[latin1]{inputenc}

\usepackage{fancyhdr}
\usepackage{times}

\usepackage{enumerate}
\usepackage{titlesec}
\usepackage{mathrsfs}
\usepackage{stmaryrd}

\def\R{\mathbb{R}}

\newcommand{\Ek}{E(\kappa,\tau)}

\newcommand{\ben}{\begin{enumerate}}
\newcommand{\bit}{\begin{itemize}}
\newcommand{\een}{\end{enumerate}}
\newcommand{\eit}{\end{itemize}}

\newcommand{\ed}{\end{document}}

\def\cR{\mathcal{R}}

\def\cW{\mathcal{W}}

\let\8=\infty \let\0=\emptyset

\let\landa=\lambda
\let\alfa=\alpha

\let\parc=\partial

\def\ep{\varepsilon}

\def\landa{\lambda}

\def\S{\Sigma}

\def\cte.{\mathop{\rm cte.}\nolimits}

\def\E{\mathbb{E}}

\def\R{\mathbb{R}}

\def\H{\mathbb{H}}
\def\S{\mathbb{S}}

\def\SS{\Sigma}

\def\Ek{\mathbb{E}^3 (\kappa,\tau)}

\def\sr{\mathbb{S}^2\times\mathbb{R}}

\newfont{\bb}{msbm10 at 12pt}

\titleformat{\section}
{\filcenter\bfseries\large} {\thesection{.}}{0.2cm}{}
\titleformat{\subsection}[runin]
{\bfseries} {\thesubsection{.}}{0.15cm}{}[.]
\titleformat{\subsubsection}[runin]
{\em}{\thesubsubsection{.}}{0.15cm}{}[.]

\usepackage[up,bf]{caption}

\newtheorem{theorem}{Theorem}[section]
\newtheorem{lemma}[theorem]{Lemma}

\newtheorem{remark}[theorem]{Remark}

\newtheorem{definition}[theorem]{Definition}

\theoremstyle{definition}

\numberwithin{equation}{section}
\numberwithin{figure}{section}

\usepackage{color}

\begin{document}

\fancyhead[LO]{Elliptic Weingarten surfaces in $\sr$}
\fancyhead[RE]{Isabel Fern\'andez}
\fancyhead[RO,LE]{\thepage}

\thispagestyle{empty}

\begin{center}
{\bf \LARGE Rotational elliptic Weingarten surfaces in $\S^2\times\R$\\[0,2cm] and the Hopf problem}
\vspace*{5mm}

\hspace{0.2cm} {\Large Isabel Fern\'andez}
\end{center}

%



\footnote[0]{
\mbox{} \hspace{0.25cm} \emph{Keywords}: Weingarten surfaces, fully nonlinear elliptic equations, phase space analysis, isolated singularities, rotational surfaces, Hopf theorem, product spaces, homogeneous spaces.}



\vspace*{7mm}

\begin{quote}
{\small
\noindent {\bf Abstract}\hspace*{0.1cm}
We prove that, up to congruence,  there exists only one immersed sphere satisfying a given uniformly elliptic Weingarten  equation in $\sr$, and it is a rotational surface.  This is obtained by showing that rotational uniformly elliptic Weingarten surfaces in $\sr$ have bounded second fundamental form together with a Hopf type result by J. A. G\'alvez and P.  Mira.  
\vspace*{0.1cm}

}
\end{quote}

%
%
%
%
%
%
%
%


\section{Introduction}

An immersed oriented surface $\Sigma$ in a Riemannian $3$-manifold $M$ is called a \emph{Weingarten surface} if its principal curvatures $\kappa_1,\kappa_2$ satisfy a smooth relation 
 \begin{equation}\label{eq:w}
 W(\kappa_1,\kappa_2)=0,
 \end{equation} 
for some $W\in C^1(\R^2)$ symmetric (i.e., $W(x,y)=W(y,x)$). The relation \eqref{eq:w} defines a fully nonlinear PDE when we view $\Sigma$ as a local graph, and we say that the Weingarten equation \eqref{eq:w} is \emph{elliptic} if this equation is elliptic.  In terms of the function $W$, this means that 
\begin{equation}\label{eq:welliptic}
\frac{\parc W}{\parc k_1}\frac{\parc W}{\parc k_2}>0 \quad \mbox{on}\; \,  W^{-1}(0). 
\end{equation}

This class of surfaces is often also referred to as {\em special Weingarten} surfaces.  The Weingarten equation \eqref{eq:w} is said to be {\em uniformly elliptic} if there exists positive constants $M_1,M_2\in\R$ such that 
\begin{equation}\label{eq:welliptic_unif}
0<M_1<\frac{\parc W}{\parc k_1}\frac{\parc W}{\parc k_2}<M_2 \quad \mbox{on}\; \,  W^{-1}(0). 
\end{equation}

The most well-known elliptic Weingarten surfaces are constant mean curvature (CMC) surfaces, and minimal surfaces in particular, for which the underlying PDE is quasilinear. In this sense,  elliptic Weingarten surfaces represent the natural \emph{fully nonlinear} extension of CMC surface theory, and an interesting problem in this setting is to explore which global results of CMC surface theory extend to the general case of elliptic Weingarten surfaces.


Regarding the ambient manifold, the most studied case is, of course,  when $M$ is a space of constant curvature.  In particular, the celebrated Hopf theorem was extended to Weingarten surfaces by Hopf \cite{H}, Chern \cite{Ch},  Hartman and Wintner \cite{HW} (see also \cite{B})  and G\'alvez and Mira \cite{GM1}, showing that any immersed elliptic Weingarten {\em sphere} (i.e., compact surface with genus $0$)  in $\R^3, \S^3$ or $\H^3$ is a  totally umbilical sphere.

After the spaces of constant curvature,  a natural family of ambient spaces to consider are the product spaces $\sr$ and $\H^2\times\R$,  which belong to the class of  simply connected  homogeneous $3$-manifolds with a $4-$dimensional isometry group, also called $\mathbb{E}^3(\kappa,\tau)$ spaces.  The Hopf problem for CMC surfaces in this setting was solved by Abresch and Rosenberg in \cite{AR, AR2}, showing that any immersed CMC sphere in a $\mathbb{E}^3(\kappa,\tau)$ space must be rotational and congruent to a unique rotational surface. 

The extension of this result for elliptic Weingarten surfaces was studied by G\'alvez and Mira in  \cite{GM2}, where they show that the Hopf problem can be solved in the affirmative (i.e., there exists a unique elliptic Weingarten immersed sphere in $\E^3(\kappa,\tau)$ and it must be rotational) if a certain rotational surface (the {\em canonical example}) has bounded second fundamental form (see Section \ref{sec:fases} for a more precise statement).  As the authors prove, this  is the case  when the ambient space is $\H^2\times\R$, but not in $\sr$,  since  the canonical example in $\sr$ for an arbitrary elliptic Weingarten  equation can have singularities, as showed in Example 8.6 in  \cite{GM2}.

Let us point out here that a remarkable difference between (general) elliptic Weingarten surfaces and CMC surfaces is that the first  ones admit singularities,  even in the euclidean case, which is not possible in the case of CMC surfaces.  In a recent work, Mira and the author \cite{FM} classified all the rotational elliptic Weingarten surfaces in $\R^3$, extending the previous classification  by Sa Earp and Toubiana \cite{ST,ST2} of complete surfaces to the singular case.  When the Weingarten equation is of {\em minimal type} (i.e., its umbilical constant is zero, see Section \ref{sec:prelim})  complete rotational elliptic surfaces in $\sr$ and $\mathbb{H}^2\times\R$ were classified in \cite{MR}, where some examples with singularities are also exhibited.



In the present paper, we will show that  uniformly elliptic rotational Weingarten surfaces in $\sr$  do not admit singularities and have bounded second fundamental form (see Theorem \ref{th:unif}).  As a byproduct of this fact together with \cite{GM2} we solve the Hopf problem in $\sr$ for  uniformly elliptic Weingarten surfaces (see Theorem \ref{th:hopf}): 

{\bf Theorem:} {\em 
There exists only one (up to congruence) immersed sphere satisfying a given uniformly elliptic Weingarten equation in $\sr$.  Moreover, this unique surface is rotational. }

The paper is organized as follows: in Section \ref{sec:prelim} we rewrite the Weingarten equation \eqref{eq:w} in a more convenient way for our purposes.  Section \ref{sec:fases} is devoted to study of the phase space associated to the elliptic Weingarten equation for rotational  surfaces in $\sr$, following the spirit of \cite{FM}.  Finally,  in Section \ref{sec:unif} we prove  Theorems  \ref{th:unif} and \ref{th:hopf} .

The author is grateful to Pablo Mira for many valuable discussions during the preparation of this paper.


\section{The elliptic Weingarten equation}\label{sec:prelim}


Let $\Sigma$ be an oriented surface in a Riemannian 3-manifold whose principal curvatures $\kappa_1\geq \kappa_2$ are related by an elliptic Weingarten equation \eqref{eq:w},  with $W\in C^1(\R^2)$ symmetric and  satisfying \eqref{eq:welliptic}.

Each connected component of $W^{-1}(0)$ gives rise to a different elliptic theory, see \cite{FGM} for a more detailed discussion. By \eqref{eq:welliptic}, any such component of $W^{-1}(0)$ can be rewritten as  a proper curve in $\R^2$ given by a graph
\begin{equation}\label{eq:weig}
\kappa_2=g(\kappa_1), \hspace{1cm} g'<0,
\end{equation}
where $g$ is $C^1$ in some interval of $\R$. In particular,  there exists a unique value $\alfa\in \R$ (which we call the \emph{umbilical constant} of the equation) such that $g(\alfa)=\alfa$. 

Recall that we are assuming $\kappa_1\geq \kappa_2$.  Thus,  $g$ is a monotonic bijection from an interval $I_g\subset [\alpha,\8)$, with $\alpha\in I_g$, to $J_g:=g(I_g)\subset (-\8,\alfa]$. By the monotonicity and properness of $g$, there are two possibilities for the intervals $I_g$ and $J_g$:

\begin{enumerate}
\item
$I_g=[\alfa,\8)$. In this case, $J_g=(b,\alfa]$, where $b  \in (-\8,\alfa)\cup \{-\8\}$ is given by $b=g(\8) = \lim_ {x\to \8} \, g(x),$ 
\item
$I_g=[\alfa,b)$, $b\in\R$. In this case ${\rm lim}_{x\to b^-} \, g(x)=-\8$ and $J_g=(-\8,\alfa]$.
\end{enumerate}

A change in the orientation of $\Sigma$ produces a surface satisfying a different elliptic Weingarten equation \eqref{eq:weig}, according to the following correspondence:  
\begin{equation}\label{eq:cambiorienta}
(\kappa_1,\kappa_2,g(x),I_g,J_g,\alfa)\mapsto (-\kappa_2,-\kappa_1,-g^{-1}(-x), -J_g,-I_g,-\alfa).
\end{equation}
In particular, up to a change of orientation on the surface $\Sigma$, we can assume that the Weingarten equation \eqref{eq:weig} on $\Sigma$ satisfies
\begin{equation}\label{eq:weigorient}
I_g=[\alfa,\8). 
\end{equation}
Let us point out that when $I_g=[\alpha,\8)$ and $J_g=(-\infty,\alpha]$, condition \eqref{eq:weigorient} is preserved by a change of orientation. However, for $I_g=[\alpha,b)$ or $J_g=(b,\alpha]$, with $b\in\R$, the choice \eqref{eq:weigorient} actually fixes an orientation on the surface. 	

Taking into account the previous discussion, throughout this paper we regard  elliptic Weingarten surfaces as  follows:

\begin{definition}\label{def:wg}
An \emph{elliptic Weingarten surface} is an oriented surface $\Sigma$ immersed in a Riemannian 3-manifold whose principal curvatures $\kappa_1\geq \kappa_2$ satisfy at every point the relation \eqref{eq:weig} for some $C^1$ map $g:I_g\to J_g$, where $I_g=[\alfa,\8)$ and $g(\alfa)=\alfa$ for some $\alfa\in \R$. We also denote $b:=g(\8)\in[-\8,\alfa)$.

We let $\cW_g$ denote the class of all (oriented) elliptic Weingarten surfaces in $\sr$ associated to a given function $g$ in these conditions.
\end{definition}

We point out here that an alternative way of defining Weingarten surfaces is by means of the condition $F(H,K_e)=0$, where $H,K_e$ are the mean and extrinsic curvatures, and $F\in\mathcal{C}^\infty(\R^2)$. In the elliptic case, on any connected component of $F^{-1}(0)$ we can write
\begin{equation}\label{eq:phi}
H=\phi(H^2-K_e),\quad \phi(t)\in C^{\8}([0,\8)),\; 4t\phi'(t)^2<1. 
\end{equation}
For example, this is the formulation used by Sa Earp and Toubiana in \cite{ST,ST2} for surfaces in $\R^3$, Morabito and Rodr\'{i}guez in \cite{MR} for surfaces in  $\sr$ and $\mathbb{H}^2\times\R$, and G\'alvez and Mira in \cite{GM2} for surfaces in the homogeneous spaces $\E^3(\kappa,\tau)$. 



\section{Rotational elliptic Weingarten surfaces in $\S^2\times\R$}\label{sec:fases}

We will regard $\S^2$ as the unit ball in $\R^3$, so $\sr$ will be seen as 
$$\sr=\{(x_1,x_2,x_3,t)\in\R^4 \; :\; x_1^2+x_2^2+x_3^2=1\}.$$

Let $\SS$ be a rotational surface in $\sr$. Up to an isometry, we can assume that $\SS$ is given by the rotation of the curve $\gamma:I\subset \R\to\sr$ 
\begin{equation}\label{eq:weigamma}
\gamma(s)=( \sin( x(s)), 0,\cos( x(s)), t(s)), \quad x(s)\in[0,\pi],\; t(s)\in\R, 
\end{equation}
around the axis $\{(0,0,1)\}\times\R$. That is, we parameterize the surface as 
$$\Sigma=\{  (\sin x(s) \cos\theta, \sin x(s) \sin \theta, \cos x(s), t(s) ) \;:\; s\in I,\theta\in [0,2\pi]\}.$$
We will also assume that the curve $\gamma$ is parameterized by arc length: $x'(s)^2 + t'(s)^2=1$. 

Observe that any rotational surface in $\sr$ is actually invariant by rotations around two axes. In this case, $\SS$ is also invariant under rotations around the {\em antipodal} axis $\{(0,0,-1)\}\times\R$. The points where $\cos x= 1$  (resp. $\cos x = -1$) correspond to the points where the surface meets the fixed rotation axis (resp. the {\em antipodal} rotation axis). 

Up to reversing the orientation of $\gamma$, we can assume that the normal vector $N$  along $\gamma(s)$ is given by 
$ N(\gamma(s))=(-t'(s)\cos x(s),0,t'(s) \sin x(s),x'(s)).$
The principal curvatures of $\SS$ are then given by: 
\begin{equation}\label{eq:curvprin}
\landa(s)=t'(s) \cot( x(s)),\qquad \mu(s)=t''(s)x'(s) - t'(s)x''(s) = \frac{t''(s)}{x'(s)} = -\frac{x''(s)}{t'(s)}.
\end{equation}

Assume now that $\SS$ is an elliptic Weingarten surface in $\cW_g$ (see Definition \ref{def:wg}). Then $\SS$ satisfies $\kappa_2=g(\kappa_1),$ where $\kappa_1\geq \kappa_2$ are the principal curvatures. Write $\SS=\SS_1\cup\SS_2$, where $\SS_i=\{p\in\SS\;:\; \kappa_i(p)=\landa(p)\}$, $i=1,2$. Then $\mu=g(\landa)$ on $\Sigma_1$ and $\landa=g(\mu)$ on $\Sigma_2$, so we can write 
\begin{equation}\label{eq:weif}
\mu=f(\landa)
\end{equation}
on $\SS$, where  $f:I_f\to I_f$, $I_f:=I_g\cup J_g$, is defined as 
\begin{equation}\label{eq:f}
f|_{I_g}=g \qquad f|_{J_g}=g^{-1}. 
\end{equation}
By virtue of \eqref{eq:weigorient}, $I_f$ is of the form $I_f=(b,\8)$, where $b=g(\8)$. 
The function $f$ is $C^1(I_f)$, and strictly decreasing, with $f \circ f = {\rm{Id}}$. 


As a consequence of \eqref{eq:curvprin} we have
\begin{equation}\label{eq:edo}
\landa'  = x'  \cot( x ) \Big(    f(\landa ) - {\landa } (1+\tan^2(x)) \Big).
\end{equation}

The above differential equation is singular when $x\in\{0,\pi\}$ and when $x=\pi/2$. The first situation happens when the generating curve $\gamma$ reaches any of the two rotation axes, $\{(0,0,1)\}\times\R$ and its {\em antipodal} axis, $\{(0,0,-1)\}\times\R$;  the second one corresponds to the case where   $\gamma$ meets the axis $\{(1,0,0)\}\times\R$. 

In particular, $(x(s),\landa(s))$ is a solution to the following nonlinear autonomous system on any open interval where $x'(s)\neq 0$ and $x(s)\in (0,\pi), $ $x(s)\neq\pi/2$:
\begin{equation}\label{eq:auto}
\left\{
\def\arraystretch{1.5}
\begin{array}{lll} 
x' & = & \ep \sqrt{1-\landa^2 \tan^2(x)}, \\ 
\landa' & = & \ep \sqrt{1-\landa^2 \tan^2(x)}  \cot( x ) \big(    f(\landa ) - \landa (1+\tan^2(x))  \big) ,
\end{array}
\right.
\end{equation}
where $\ep ={\rm sign} (x') =\pm 1$. This process can be reversed, so that any solution to \eqref{eq:auto} with $x(s)\in (0,\pi),$ $x(s)\neq \pi/2$, determines a rotational surface $\Sigma$ of the Weingarten class $\cW_g$. Thus, the orbits of \eqref{eq:auto} will be identified with the profile curves of rotational surfaces in $\cW_g$ on open sets where $x'(s)\neq 0$ and $x(s)\in (0,\pi)$, $x(s)\neq\pi/2$.

\begin{remark}\label{re:eps}
The systems \eqref{eq:auto} for $\ep = 1$ and $\ep = -1$ are actually equivalent, since they have the same orbits. Indeed, if $(x(s),\landa(s))$ is a solution for $\ep = 1$, then $(x(-s), \landa(-s))$ is a solution for $\ep=-1$.  
\end{remark}

\begin{remark}\label{re:pi2}
Notice also that if $(x(s),\landa(s))$ is a solution of \eqref{eq:auto}, then $(\pi-x(s),\landa(s))$ is also a solution. The corresponding rotational surfaces differ by an isometry of $\sr$ interchanging the rotational axes $\{(0,0,1)\}\times\R$ and $\{(0,0,-1)\}\times\R$.  
\end{remark}

\begin{figure}[htbp]
\begin{center}
\includegraphics[width=.5\textwidth]{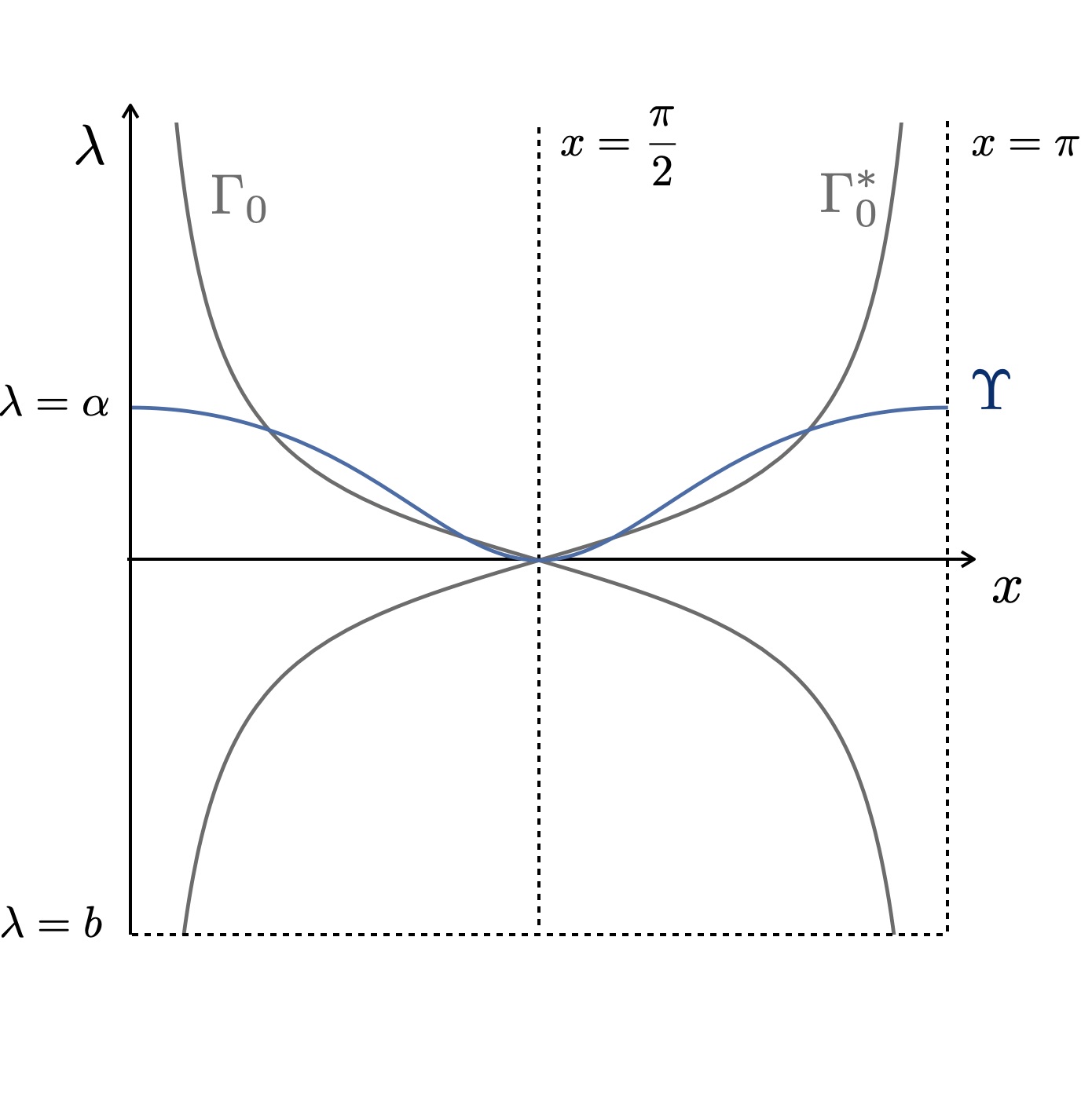}\caption{Phase space of \eqref{eq:auto} and the curve $\Upsilon$ (here $\alpha>0>b>-\8$). } \label{fig:Upsilon}
\end{center}
\end{figure}

The phase space of \eqref{eq:auto} is the set $\cR=\cR_0\cup\cR_0^\ast$, where 
$$\cR_0=\{ (x,\landa)\in \R^2:\; x\in (0,\pi/2),\; \landa >b,\; |\landa\tan(x)|<1 \}$$
$$\cR_0^\ast=\{ (x,\landa)\in \R^2:\; x\in (\pi/2,\pi),\; \landa>b, \;  |\landa\tan(x)|<1 \}.$$
We also denote by $\Gamma=\Gamma_0\cup\Gamma_0^\ast$ the  boundary curves given by 
\begin{eqnarray}\label{eq:Gamma}
 \Gamma_0 &:=& \{(x,\landa) : 0<x<\pi/2,  \; \landa>b, \; |\landa \tan(x)| = 1\}\\
 \Gamma_0^\ast &:=& \{(x,\landa) : \pi/2<x<\pi,    \; \landa>b, \; |\landa \tan(x)| = 1\}
\end{eqnarray}
(see Figure \ref{fig:Upsilon}).

\begin{remark}\label{re:simetria}
$\Gamma$ corresponds to the points where the generating curve $\gamma(s)\equiv (x(s),t(s))$ of $\Sigma$ does not intersect the  axis $\{1,0,0\}\times\R$ and has  vertical tangent vector. 
By Remark \ref{re:eps}, if $\gamma(s)$ has a point $s_0\in\R$ with $x(s_0)\neq \pi/2$ and vertical tangent vector, its associated orbit $(x(s),\landa(s))$ in \eqref{eq:auto} \emph{hits} $\Gamma$ at $s=s_0$, and then {\em bounces back}, following the same trajectory in the opposite sense, but with the sign of $\ep$ reversed. Thus,  $\gamma(s)$ extends smoothly across $s_0$, with the sign of $x'(s)$ changing at $s_0$, and  $\Sigma$ is symmetric with respect to the horizontal section passing through $\gamma(s_0)$.

Moreover, any orbit $(x(s), \landa(s))$ hits $\Gamma$ at most at two points.  In particular, $x'(s)\neq 0$ except at these points, and $(x(s),\landa(s))$ is a graph in $\cR$ over the $x$-axis.  
\end{remark}

We now describe some properties regarding the behavior of the orbits of \eqref{eq:auto} that will be useful in our study.  
In the sequel, we will denote by $\Upsilon$ the 
curve in the $(x,\landa)$-plane   given by the equation 
\begin{equation}\label{eq:upsilon}
f(\landa)={\landa} (1+\tan^2(x)), \quad x\in (0,\pi),\; \landa>b,
\end{equation} 
where $f$ is as in \eqref{eq:f}. 
Since $f$ is strictly decreasing with $f(\alpha)=\alpha$ and $f\circ f=Id$, it follows that
if the umbilical constant $\alpha$ is positive, $\Upsilon$ lies in the band $\landa\in[0,\alpha]$, whereas for $\alpha<0$ the curve $\Upsilon$ lies in $\landa\in[\alpha,0]$, and if $\alpha=0$ then $\Upsilon$ is given by  $\{\landa=0\}$.  Moreover,  on $\Upsilon$ we have 
$$ \big(   f'(\landa) - 1 - \tan^2(x) \big)  {d\landa} =\big (  2\landa \tan(x) (1 + \tan^2(x))   \big) dx,$$
and therefore 
$\Upsilon$ is the graph of a function $h(x)$ defined on $(0,\pi)$,  symmetric with respect to $x=\pi/2$,  with $h(0)=h(\pi)=\alpha$.   
When $\alpha > 0$ (resp. $\alpha<0$), $h(x)$ is strictly decreasing (resp. increasing) on $(0,\pi/2)$ and strictly increasing (resp. decreasing) on $(\pi/2,\pi)$. 
At $x=\pi/2$,  we have $h(\pi/2)=\max\{0,b\}$.  Indeed, by \eqref{eq:upsilon} we have $f(h(x))/h(x) \to \8$ when $x\to\pi/2$. Thus, either $h(x)\to 0$ and $b<0$ (observe that $f(\lambda)$ is only defined for $\lambda>b$),  or $h(x)\to b$ (and therefore $f(h(x))\to\8$) and $b>0$.  
Finally,  notice that the restriction of $\Upsilon$ to the phase space $\cR$ may or may not be connected (see Figure \ref{fig:Upsilon}). 

The curve $\Upsilon$ determines the monotonicity regions of the phase space $\cR$, as stated in the following lemma. 

\begin{lemma}\label{lem:monot}
Any orbit of \eqref{eq:auto}  is a graph  in $\cR$ over the $x$-axis of a function that is strictly decreasing (resp. increasing) whenever the orbit lies above (resp. below)  $\Upsilon\cap\cR_0$, and is strictly increasing (reap. decreasing) whenever it lies above (resp. below) $\Upsilon\cap\cR_0^\ast$.  In particular, any orbit in $\cR_0$ (resp. in $\cR^\ast_0$) intersects transversely $\Upsilon$ at most once,  unless $\alpha=0$ and the orbit is given by $\landa\equiv 0$, that coincides with $\Upsilon$. 
\end{lemma}

\begin{proof}
The fact that $(x(s),\landa(s))$ is a graph over the $x$-axis was discussed in Remark \ref{re:simetria}. The monotonic character of this graph follows from a careful but direct analysis  of \eqref{eq:edo}.  
In particular, this monotonicity behavior shows that if $(x(s),\landa(s))$ intersects $\Upsilon$ transversely at some point $(x_0,\landa_0)\in\cR_0$,  the orbit must lie below $\Upsilon$ for $x\to x_0^-$, and above $\Upsilon$ for $x\to x_0^+$ in the case $\alpha>0$ (recall that, for  $\alpha>0$ and $x\in (0,\pi/2)$,  $\Upsilon$ is given by the graph of a strictly decreasing function); whereas if $\alpha<0$  the orbit must lie above $\Upsilon$ for $x\to x_0^-$ and below $\Upsilon$ for $x\to x_0^+$ . In particular, this shows that $(x(s),\landa(s))$  intersects transversely the curve $\Upsilon$ in $\cR_0$ at most once.  The case of $\cR_0^\ast$ is similar. 

\end{proof}


\begin{lemma}\label{lem:b}
Let $(x(s),\landa(s))$ be an orbit of \eqref{eq:auto} where $f$ is as in \eqref{eq:f} with $b\neq -\8$. Assume that  the orbit meets $\partial\cR$ at a point $(x_0,b)$. Then $x_0\notin \{ 0,\pi \}$. 
\end{lemma}

\begin{proof}
Arguing  by contradiction, assume $(x(s),\landa(s))$ approaches to the point $(0,b)$ as $s\to s_0$ (the case $(\pi,b)$ is analogous, see Remark \ref{re:pi2}). Since $b\in\R$,  $\landa(s)=t'(s)\cot(x(s))$ is bounded as $s\to s_0$ and, in particular, $t'(s)\to 0$ as $s \to s_0$ . The profile curve $\gamma(s)\equiv  (x(s),t(s))$ of the rotational elliptic Weingarten surface $\Sigma$ associated to the orbit can then be reparameterized, for $x$ close to $0$, as the graph of a function $h(x)$ with $h'(x)\to 0$ as $x\to 0$. The principal curvatures of $\Sigma$ are  given as 
$$\landa=\pm h'(x)(1+h'(x)^2)^{-1/2}\cot(x), \qquad \mu=f(\landa)=(1+h'(x)^2)^{-1}h''(x),$$
where the $\pm$ sign is given by the sign of $x'(s)$ as $s\to s_0$.
In particular,  $\lim_{x\to 0}h''(x)= \lim_{\landa\to b} f(\landa)= \8,$ which contradicts that $\lim_{x\to 0} \frac{h'(x)}{\sin(x)}=\pm b$. 

\end{proof}


\subsection{The canonical example}\label{sub:canonical}

In \cite{GM2} it is proved that, given an elliptic Weingarten class $\cW_g$  of surfaces in $\sr$ (and more generally, in any $\Ek$ space), there exists a unique extensible, regular rotational surface $\Sigma_0$  meeting orthogonally the rotation axis. This surface will be referred to as the {\em canonical example} of the Weingarten class $\cW_g$.  As pointed out in \cite[Example 8.6]{GM2},  the canonical example in  $\sr$ can have singularities (this situation does not happen in $\R^3$,  where the canonical example corresponds to a totally umbilical sphere/plane of principal curvatures equal to $\alpha$).  

When the canonical example has bounded second fundamental form, the following result regarding the classification of Weingarten (topological) spheres in the $\Ek$ spaces is obtained: 

\begin{theorem} \cite[Theorem 1.6]{GM2} \label{th:GM}
 Let $\cW_g$ be the class of elliptic Weingarten surfaces in $\sr$ given by \eqref{eq:weig}.  Assume that the canonical example $\Sigma_0$ in $\cW_g$ has bounded second fundamental form. Then, any immersed topological sphere $\Sigma$ in $\cW_g$ is a rotational sphere.  More specifically, if $\Sigma_0$ is compact, then $\Sigma$ is congruent to $\Sigma_0$.  If  $\Sigma_0$ is not compact, then $\Sigma$ does not exist.
\end{theorem}

Observe that, since any point lying on the rotational axis of a (regular) rotational surface in $\sr$  must be umbilical,  the  associated orbit of the canonical example $\Sigma_0$ reaches the boundary of the phase space at $x=0,\landa=\alfa$. Moreover, if $\alpha=0$ then $\Sigma_0$ is the slice $\S^2\times\{0\}$ corresponding to the orbit $\landa\equiv 0$, whereas for $\alfa\neq 0$, the corresponding orbit is above the curve $\Upsilon$ given by \eqref{eq:upsilon} (see  Lemma \ref{lem:monot}). 



\subsection{Singularities for rotational surfaces in $\cW_g$}


Let $\Sigma$ be a rotational surface in $\cW_g$ given by the rotation of an arc-length parameterized curve $\gamma(s)$ as in \eqref{eq:weigamma} and defined in some interval $I\subset \R$. Let $a\in\R$ be one of the endpoints of $I$ and assume that $\Sigma$ has a singularity at $s=a$. That is, $\Sigma$ cannot be extended to a regular surface on $I\cup\{a\}$. 

The following lemma characterizes the behavior of the singularities in terms of the corresponding orbit for the system \eqref{eq:auto}. 

\begin{lemma}\label{lem:sing}
Let $\Sigma$ be a rotational surface in $\cW_g$ and 
$(x(s),\landa(s))$ its associated orbit in $\cR$ for \eqref{eq:auto}. 
Then, $\Sigma$ has a singularity 
at $s= a$
if and only if one of the following two conditions holds:
\begin{enumerate}[i)]
\item 
$(x(s),\landa(s))\to (x_0,\pm\8)$ as $s\to a$, with $x_0 \in \{0,\pi\}$  (the case $\landa(s)\to-\8$ can  only  happen if $b=g(\8)=-\8$), or 
\item  
$(x(s),\landa(s))\to (x_0,b)$ as $s\to a$, with $x_0\in (0,\pi)$ and $b\neq -\8$. 
\end{enumerate}

The first case corresponds to an isolated singularity created as the surface touches one of its two rotational axes. In the second one, the surface is singular along a compact curve. 

\end{lemma}

\begin{proof}
If the associated orbit to $\Sigma$ satisfies $i)$ or $ii)$ above it is clear that the second fundamental form of $\Sigma$ is not bounded as $s\to a\in\R$, and so the surface has a singularity.

Conversely,  assume that  $\Sigma$ has a singularity at $s= a\in\R$, and therefore its associated orbit $(x(s),\landa(s))$  approaches the boundary of the phase space $\cR$. Moreover, since the orbit is a graph over the $x$-axis,  $x(s)$ has a well-defined limit $x_0\in [0,\pi]$ as $s\to a$. As discussed in Remark \ref{re:simetria}, $\Sigma$ can be smoothly extended if $(x(s),\landa(s))$ approaches  $\Gamma$. Thus, have the following possibilities: 

\begin{enumerate}[a)]

\item $x(s)$ converges to $x_0= \pi/2$ as $s\to a$. By the first equation in  \eqref{eq:curvprin} we have that $0\in \overline{I_f}$ and $\landa(s) \to 0$ as $s\to a$. If $0\in I_f$ then  $\mu(s)\to f(0)\in\R$ and by the second equation in \eqref{eq:curvprin} the profile curve $\gamma(s)\equiv (x(s),t(s))$ of $\Sigma$ would have bounded second derivatives as $s\to a$. In particular, the mean value theorem implies that  $\gamma'(s)$ and $\gamma(s)$ have well defined limits $(x'_0,t'_0)$ and $(x_0,t_0)$ as $s\to a$, and $\gamma(s)$ can be extended across $s=a$ as the (unique) solution to the Cauchy problem:
\begin{eqnarray*}
x''(s) & = & -t'(s) \, f(t'(s) \cot x(s))\\
t''(s) & = & x'(s) \,  f(t'(s) \cot x(s) ) 
\end{eqnarray*}
with initial conditions $(x(a),t(a))=(x_0,t_0)$, $(x'(a),t'(a))=(x'_0,t'_0)$, contradicting that $\Sigma$ has a singularity at $s=a$. Therefore, $\landa(s)\to 0=b$, as we wanted to prove.  

\item The orbit $(x(s),\landa(s))$ approaches $\{x=0\}\cup\{x=\pi\}$. In this case, it remains to prove that $\landa(s)\to\pm\8$. By $\eqref{eq:edo}$ and Lemma \ref{lem:monot},   $\landa(s)$ is monotonic as $x(s)\to x_0\in\{0,\pi\}$ and so it has a well defined limit $\landa_0$ as $s\to a$.  The case $\landa_0=b\neq -\8$ is impossible by Lemma \ref{lem:b}.  If  $\landa_0\neq \pm\8$, the mean curvature of $\Sigma$ would extend continuously to the (isolated) singularity with the value $1/2(\landa_0+f(\landa_0))$, which is impossible by \cite{LR}.  Thus, $\landa_0=\pm\8$ and we are done. 

\item The orbit $(x(s),\landa(s))$ approaches  $\{\landa=b\}$ in $\cR$ (when  $b\neq -\8$). In this case,  Lemma \ref{lem:b} shows that $x_0\notin \{0,\pi\}$, finishing the proof. 
\end{enumerate}



\end{proof}


\section{The uniformly elliptic case} \label{sec:unif}

Let $\cW_g$ be the Weingarten class given by \eqref{eq:weig}.  
The uniform ellipticity condition \eqref{eq:welliptic_unif} can be rephrased in terms of $g$  as
\begin{equation}\label{eq:unif}
\kappa_2 = g(\kappa_1) ,\qquad \Lambda_1 < g'< \Lambda_2<0
\end{equation} 
for some constants $\Lambda_1,\Lambda_2<0$ (see e.g. \cite{FGM,FM}).  In particular, $b=g(\8)=-\8$.  Our previous analysis in Section \ref{sec:fases} can be used to prove that rotational, uniformly elliptic Weingarten surfaces do not present singularities:

\begin{theorem}\label{th:unif}
Let  $\Sigma$ be a rotational uniformly elliptic Weingarten surface in $\sr$. Then,  $\Sigma$ has bounded second fundamental form and it cannot have singularities. 
\end{theorem}

\begin{proof} 
Let $(x(s), \landa(s))$ in $\cR$ be the solution of \eqref{eq:auto} associated to the surface $\SS$. By Lemma \ref{lem:monot}, the orbit is a graph over the $x$-axis $(x,\landa(x))$ satisfying \eqref{eq:edo}. 

Assume that the second fundamental form of $\SS$ is not bounded. Taking into account that in the uniformly elliptic case we have $b=g(\8)=-\8$, this can only occur if $\landa\to \pm\8$ and therefore $x\to 0$ or $x\to\pi$. More specifically, we have $|\tan(x)\landa(x)|<1$  as $x\to 0$ or $x\to \pi$, with $|\landa(x)|\to\8$. 
By Remark \ref{re:pi2}, it suffices to check the case where $x\to 0$. We will also assume that $\landa(x)\to\8$ as $x\to 0$ (the case $\landa(x)\to -\8$ is analogous). 

As a consequence of \eqref{eq:unif} there exist $A_1,A_2\in\R$, and $M_1,M_2<0$ such that 
\begin{equation}\label{eq:lineal}
f_1(r) < f(r) < f_2(r),
\end{equation}
where $f_i(r)=M_i r + A_i$, $i=1,2$. Moreover, $A_i=\alpha (1-M_i)$.


Let $(x_0,\landa_0)\in\cR_0$ be any point in the orbit $(x,\landa(x))$ and label $(x, \landa_i(x) )$ as the corresponding orbits of \eqref{eq:auto} for the Weingarten relations given by $f_i$, $i=1,2$ and passing through $(x_0,\landa_0)$. Without loss of generality, we can assume that the fixed point $(x_0,\landa_0)$ is not cointained in the orbit associated to the canonical example for the Weingarten relation given by  $f_2$. 

Set $y=\tan x$. Then \eqref{eq:edo} leads to 
\begin{equation}\label{eq:edo2}
\landa'(y)=\frac{1}{y} \Big(\frac{f(\landa)}{1+y^2} - \landa \Big). 
\end{equation}

Equation \eqref{eq:edo2} can be explicitly solved in the linear case, $f_i(r)=M_i r + A_i$, for which we have  
\begin{equation}\label{eq:explicit}
\landa_i(y)=\frac{1}{2} (1+y^2)^{-M_i/2} \Big(    \left(-y^2\right)^{\frac{M_i-1}{2}}A_i B_i(-y^2) + 2 C_i y^{M_i-1}            \Big)
\end{equation}
for suitable $C_i\in\R$, $i=1,2$.  Here $B_i(r)$ denotes the {\em incomplete beta function}: 
$$B_i(r):=B(r;a_i,b_i) = \int_0^r s^{a_i-1} (1-s)^{b_i-1} ds$$
for $a_i=(1-M_i)/2,$ $b_i=M_i/2$. We remark here that $B_i(r)$ satisfies 
\begin{equation}\label{eq:beta}
\lim_{r\to 0} r^{\frac{M_i-1}{2}} B_i(r)= \frac{2}{1-M_i}
\end{equation}
for $i=1,2$, which in particular implies $C_2\neq 0$ (otherwise $\landa_2(0)=\alpha$ and, taking into account Lemma \ref{lem:sing}, this would contradict the uniqueness of the canonical example). 


From  \eqref{eq:lineal} and taking into account \eqref{eq:edo2}, and the corresponding analogous equations for $\lambda_i'(y)$, $i=1,2$,  it follows that 
$$\landa_1'(y) < \landa'(y) < \landa_2'(y)$$ 
for $y>0$. Since $\landa_i(y_0)=\landa(y_0)=\landa_0$, $i=1,2$, then
$$\landa_2(y)\leq \landa(y)\leq \landa_1(y)$$
 for any $0<y\leq y_0$. In particular, by \eqref{eq:explicit} and \eqref{eq:beta}
 it follows from the above inequality that $|y\landa(y)|\to\8$ as $y\to 0^+$, which gives that $\SS$ cannot have unbounded second fundamental form as $x\to 0$.  
 
Finally, by Lemma \ref{lem:sing} it follows that at a singularity, either $\landa(s)\to\pm \infty$ or $\landa(s)\to b$ (and so, $\mu(s)=f(\landa(s))\to \infty$). In particular, as the second fundamental form of $\Sigma$ is bounded, the surface does not have singularities. 
\end{proof}

As a particular case of the above theorem, if $g$ is uniformly elliptic, the {\em canonical example} in $\cW_g$ has bounded second fundamental form. Thus, as a consequence of \cite{GM2} (see  Theorem \ref{th:GM}) we have the following Hopf-type theorem for uniformly elliptic  Weingarten surfaces in $\sr$:

\begin{theorem}\label{th:hopf}
For any given class $\cW_g$ of uniformly elliptic Weingarten surfaces
there exists only one (up to congruence) immersed sphere in $\sr$  in $\cW_g$.   Moreover, this unique surface is rotational. 
\end{theorem}

\begin{proof} 
By Theorems \ref{th:GM} and \ref{th:unif}, it only remains to check that the canonical example $\Sigma_0$ in $\cW_g$ is compact. Let $(x,\landa(x))$ be the associated orbit. Since $\Sigma_0$ has bounded second fundamental form, the orbit cannot approach $\landa\to\pm \infty$. Thus either $(x,\landa(x))\to (\pi,\alpha)$ or $(x,\landa(x))\to (x_0,\landa_0)\in \Gamma$. In the first case $\Sigma_0$ is a graph over the whole $\S^2$ section in $\sr$, and therefore is compact. In the second one, by Remark \ref{re:simetria} it suffices to prove that $(x(s),\landa(s))\to (x_0,\landa_0)$ as $s\to s_0$ with $s_0\in\R$. 

Assume, arguing  by contradiction, that $s_0=\8$. Then $t'(s)^2\to 1$ and $x''(s)=-t'(s)f(\landa(s)) \to \pm f(\landa_0)$. On the other hand, since $x(s)\to x_0$ as $s\to \8$,  by the mean value theorem  $x'(s)\to 0$ and subsequently $x''(s)\to 0$. In particular,  $f(\landa_0)=0$. By the monotonicity of $f$, it is clear that $|\alfa|<|\landa_0|$, which contradicts the behavior of the orbits described in Lemma \ref{lem:monot} and finishes the proof.  
\end{proof}



\def\refname{References}
%


\vskip 0.2cm

\noindent Isabel Fern\'andez

\noindent Departamento de Matem\'atica Aplicada I,\\ Instituto de Matem\'aticas IMUS \\ Universidad de Sevilla (Spain).

\noindent  e-mail: {\tt isafer@us.es}

\noindent This research has been financially supported by Project PID2020-118137GB-I00 funded by MCIN/AEI /10.13039/501100011033.

\end{document}